\documentclass[10pt]{amsart}
\usepackage{graphicx}
\usepackage{bbm}

\theoremstyle{plain}

\newtheorem*{theorem*}{Theorem}

\newtheorem*{conjecture*}{Conjecture}

\newcommand{\bbC}{\mathbb{C}}
\newcommand{\bbN}{\mathbb{N}}
\newcommand{\bbQ}{\mathbb{Q}}
\newcommand{\bbR}{\mathbb{R}}
\newcommand{\bbS}{\mathbb{S}}

\newcommand{\CPone}{\bbC\mathrm{P}^1}
\DeclareMathOperator{\diag}{diag}
\newcommand{\ol}[1]{\overline{#1}}
\newcommand{\half}{\tfrac{1}{2}}
\DeclareMathOperator{\tr}{tr}
\DeclareMathOperator{\Imag}{Im}
\DeclareMathOperator{\id}{\mathbbm{1}}

\newcommand{\MatrixGroup}[1]{{\mathrm{#1}}}
\newcommand{\matsl}[2]{\MatrixGroup{sl}_{#1}{#2}}

\numberwithin{equation}{section}

\newcommand{\coloneq}{:=}

\newcommand{\spacecomma}{\,\,,}
\newcommand{\spaceperiod}{\,\,.}
\newcommand{\deriv}{\mathrm{d}}

\theoremstyle{plain}

\newtheorem*{experiment*}{Experiment}

\title{Exploring the space of compact symmetric CMC surfaces}

\author{Lynn Heller}
\author{Sebastian Heller}
\author{Nicholas Schmitt}

\address{Lynn Heller \\
  Institut f\"ur Mathematik\\
  Universit{\"a}t T\"ubingen\\ Auf der Morgenstelle
  10\\ 72076 T\"uubingen\\ Germany}
\email{lynn-jing.heller@uni-tuebingen.de}

\address{Sebastian Heller \\
  Institut f\"ur Mathematik\\
  Universit{\"a}t T\"ubingen\\ Auf der Morgenstelle
  10\\ 72076 T\"uubingen\\ Germany}
\email{heller@mathematik.uni-tuebingen.de}

\address{Nicholas Schmitt \\
  Institut f\"ur Mathematik\\
  Universit{\"a}t T\"ubingen\\ Auf der Morgenstelle
  10\\ 72076 T\"ubingen\\ Germany}
\email{nschmitt@mathematik.uni-tuebingen.de}


\thanks{The first author is supported by the European Social Fund,
  by the Ministry of Science, Research and the Arts Baden-W\"urtemberg
  and by the Baden-W\"urtemberg Foundation;
  the other authors are supported by the DFG through the project HE 6829/1-1.}

\begin{document}

\begin{abstract}
  We map out the moduli space of Lawson symmetric
  constant mean curvature surfaces in the 3-sphere of genus $g>1$
  by flowing numerically from
  Delaunay tori with even lobe count
  via the generalized Whitham flow.
\end{abstract}

\maketitle

\subsection{Overview}

In this note we map out a portion of the moduli space of
embedded constant mean curvature (CMC) surfaces in the 3-sphere
experimentally
by a numerical implementation
of the generalized Whitham flow~\cite{Heller_Heller_Schmitt_2015a}.
This provides
numerical evidence for the existence of the flow
reaching arbitrary genus.

\begin{experiment*}
For each pair of integers $g\ge 1$ and $n\ge 0$ we construct
numerically a $1$-parameter family $\Xi_{g}^{n}$ of compact
Alexandrov embedded CMC surfaces of genus $g$ in
$\bbS^3$, with $n$ controlling the lobe count:
\begin{itemize}
\item
  the family $\Xi_{g}^{0}$ starts at the Lawson surface $\xi_{g,1}$
  and converges to a chain of two minimal spheres;
\item
  the family $\Xi_{g}^{n}$ ($n\ge 1$) converges at one end
  to a chain of $(g+1)n$ CMC spheres and at the other to a
  chain of $(g+1)n+2$ CMC spheres.
\end{itemize}
Each surface in $\Xi_{g}^{n}$ has a cyclic symmetry of order $g+1$
with four fixed points.
\end{experiment*}

These $\Xi_{g}^{n}$ families were computed numerically via the
generalized Whitham flow~\cite{Heller_Heller_Schmitt_2015a}, a
topology-breaking flow through CMC surfaces in $\bbS^3$ which starts
at CMC tori and, as indicated by numerical evidence,
reaches closed CMC surfaces of arbitrary genus.

\begin{figure}
\centering
\includegraphics[width=0.6\textwidth]{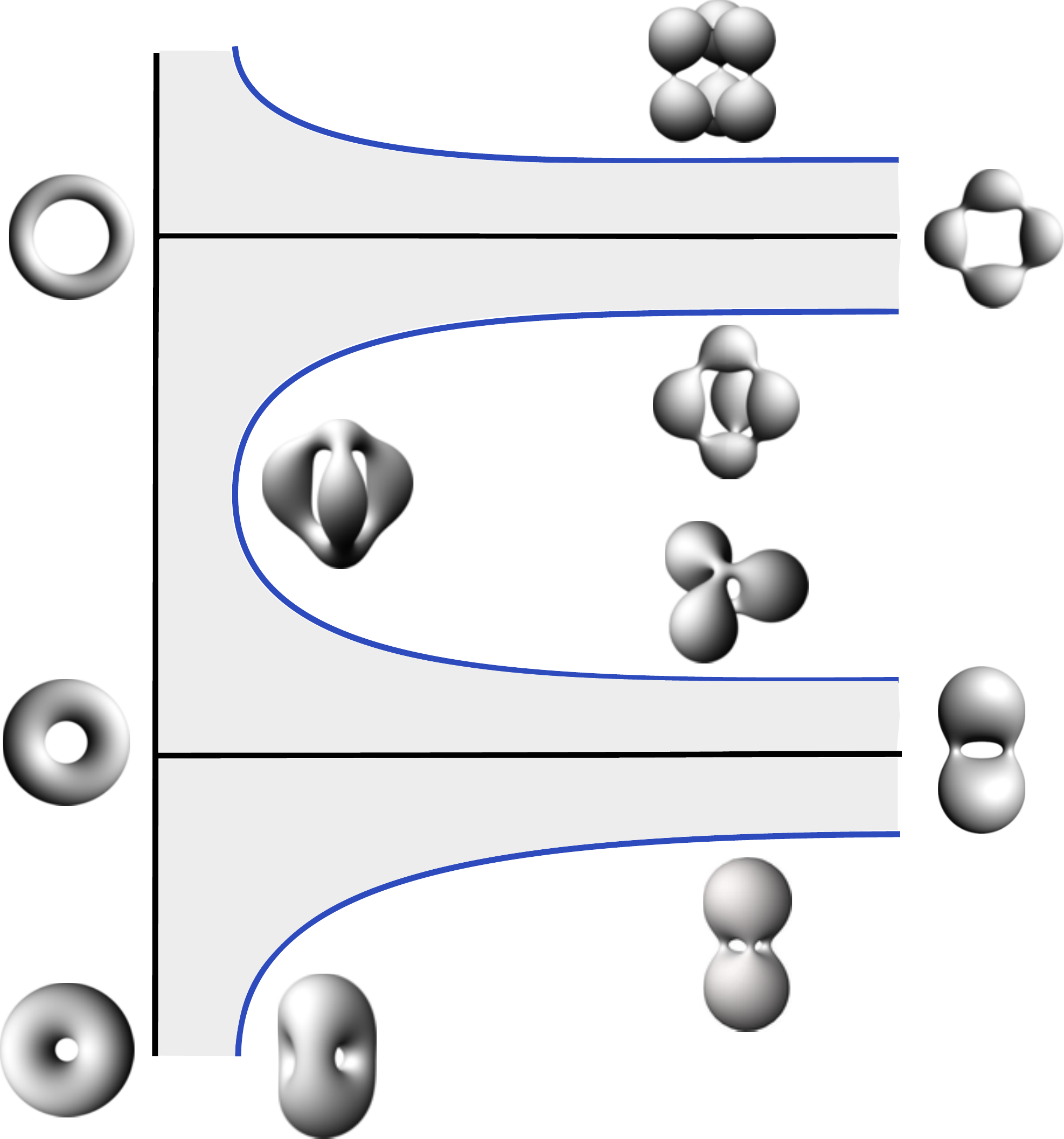}
\caption{
  \footnotesize
  The moduli space of Lawson symmetric CMC surfaces in $\bbS^3$
  arising from even-lobed Delaunay tori.
  The vertical line represents homogeneous tori of increasing mean curvature
  starting at the minimal Clifford torus at the bottom.
  The horizontal lines represent Delaunay tori
  with increasing even number of lobes (bottom to top).
  The curved lines represent the 1-parameter families
  $\Xi_{2}^{n}$ of genus $2$ surfaces.
  The flow occurs in the shaded regions;
  The families $\Xi_{g}^{n}$ of genus $g>2$
  (not shown) are obtained by continuing
  the shaded regions beyond the curved lines.
}
\label{fig:lawson-moduli}
\end{figure}

The generalized Whitham flow passes through each of the families
$\Xi_{g}^{n}$ with $n$ fixed and $g$ increasing arbitrarily, starting
at the tori $\Xi_{1}^{n}$.  To describe this initial data,
recall~\cite{Kilian_Schmidt_Schmitt_2010} the embedded CMC tori in the
3-sphere consist of the 1-parameter family of homogeneous tori of
increasing mean curvature starting at the minimal Clifford torus,
along which bifurcate 1-parameter families of $n$-lobed Delaunay
(equivariant) tori at sequencial bifurcation points $\beta_{m}$ (see
figure \ref{fig:lawson-moduli}).  The initial family $\Xi_{1}^{0}$ is
made up of the homogenous tori between the Clifford torus and
$\beta_{2}$, together with the $2$-lobed Delaunay tori.  The initial
family $\Xi_{1}^{n}$ ($n>0$) is made up of the $(2n)$-lobed Delaunay
tori, the homogeneous tori between $\beta_{2n}$ and $\beta_{2n+2}$,
and the $(2n+2)$-lobed Delaunay tori.

The topology-breaking flow is described qualitatively as follows.  The
initial torus in $\Xi_{1}^{n}$ has a cyclic symmetry of order two with
four fixed points.  The flow preserves the topology of the torus minus
two disks, formed by introducing two cuts connecting the fixed points
in pairs.  The flow retains the rotational symmetry, decreasing its
angle $\alpha$ from $\pi$ to $0$.  When $\alpha = 2\pi/(g+1)$,
$g\in\bbN$, the surface can be completed by the rotational symmetry to
a closed compact unbranched surface of genus $g$.  At other angles
$\alpha\in2\pi\bbQ$, the surface can be completed to a closed surface
branched at four points.

\begin{figure}
\centering
\includegraphics[width=0.45\textwidth]{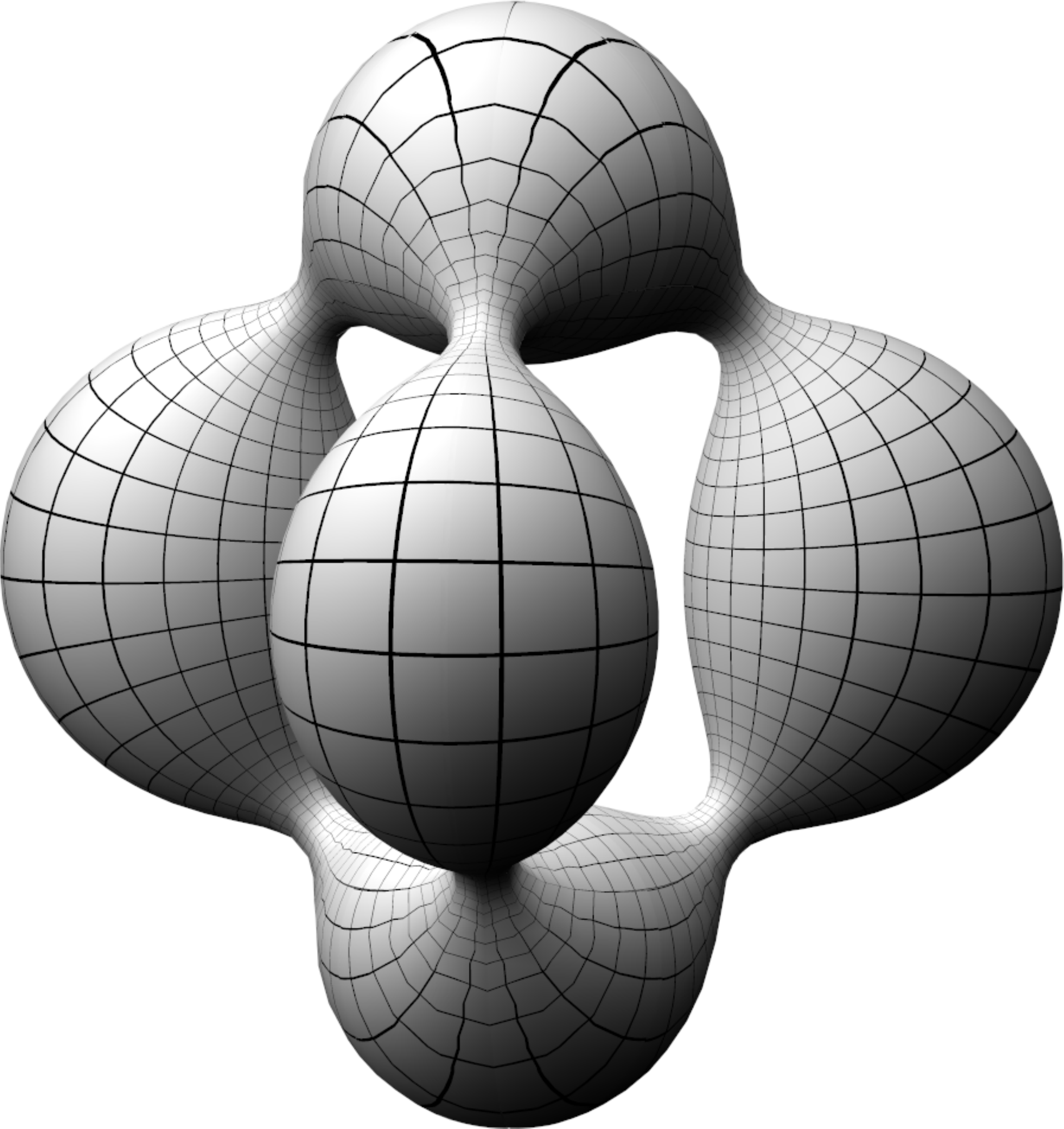}
\hspace{3ex}
\raisebox{0.03\height}{
\includegraphics[width=0.38\textwidth]{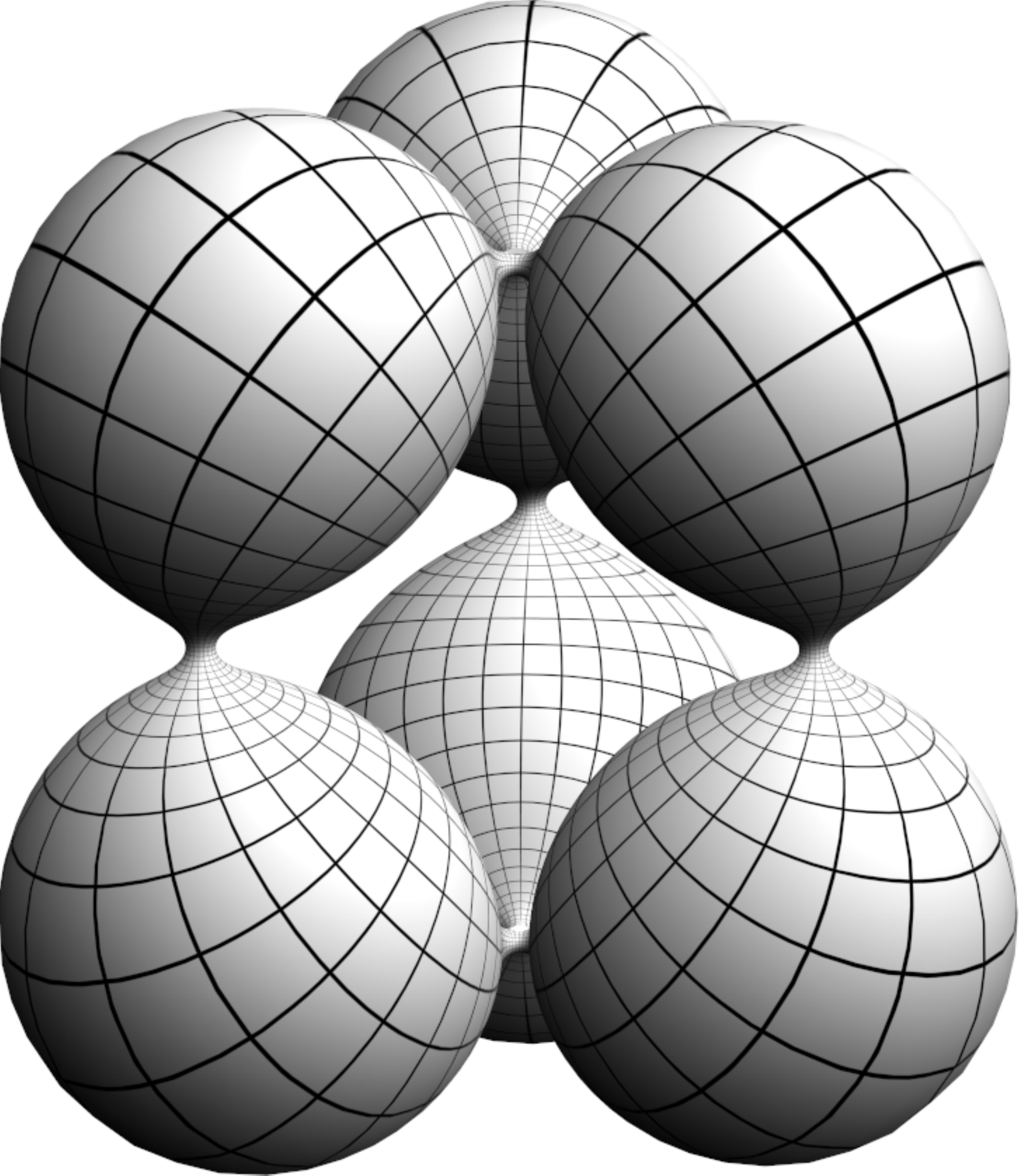}}
\caption{
  \footnotesize
  Genus 2 surface in $\Xi_{2}^{1}$ with five lobes (left).
  Genus 2 surface in $\Xi_{2}^{2}$ with six lobes (right).
}
\label{fig:surface}
\end{figure}

The generalized Whitham flow starting at a $(2n)$-lobed Delaunay torus
can flow either to $\Xi_g^{n-1}$ or $\Xi_g^{n}$. This flow direction
is determined by the choice of order two symmetry:
the symmetry with axes through necks flows to $\Xi_g^{n-1}$
while the symmetry with axes through bulges flows to $\Xi_g^{n}$
(see figure \ref{fig:surface}).

\begin{figure}
\centering
\includegraphics[width=0.4\textwidth]{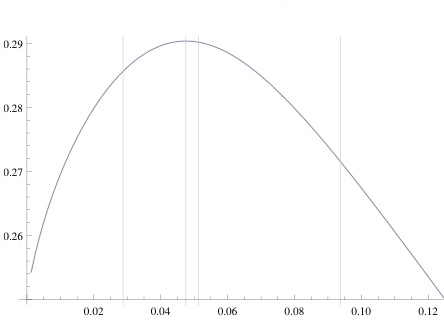}
\includegraphics[width=0.4\textwidth]{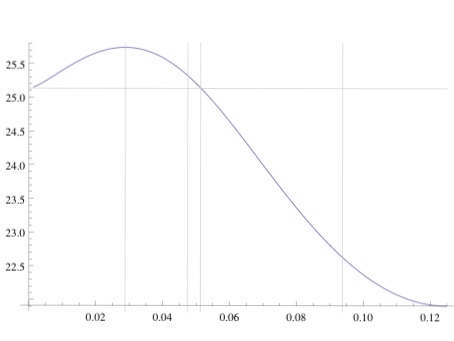}
\caption{The above plots show the mean curvature (left) and Willmore energy (right) along the  family of genus 2 CMC surfaces $\Xi_{g}^{0}$ depending on its conformal type  encoded as the poleangle $\theta\in [0, 1/8]$. For each plot the family $\Xi_{g}^{0}$ starts at the Lawson surface (right) and ends at the chain of two spheres (left): the four poles of the DPW potential (umbilics of the surface) are $\pm e^{\pm 2 \pi i \theta}$. The events from right to left are marked by vertical lines: 1. Branchpoint on unit circle, 2. Willmore energy $ = 8 \pi$, 3. Maximum mean curvature, 4. Maximum Willmore energy ($> 8 \pi$).
}
\label{fig:lawson-data}
\end{figure}

The families $\Xi_{g}^{0}$ and $\Xi_{g}^{1}$ were
first discovered in~\cite{Heller_Schmitt_2015} by numerical search. We have numerically computed the mean curvature and the Willmore energy of $\Xi_{g}^{0}$, see Figure \ref{fig:lawson-data}.

\subsection{The potential}
We construct the $\Xi_{g}^{n}$ families
numerically via the generalized Weierstrass
representation~\cite{Dorfmeister_Pedit_Wu_1998}
for CMC surface in $\bbS^3$.
The Weierstrass data consists in a $\matsl{2}{\bbC}$
loop-valued potential
$\xi$ with appropriate asymptotics in the spectral parameter $\lambda$.
The CMC immersion is obtained as
$F(\lambda_0)F^{-1}(\lambda_1)$,
where $\lambda_0,\,\lambda_1\in\bbS_\lambda^1$ are the \emph{sym points},
and $F$ is the unitary factor of the loop group Iwasawa factorization
of $\Phi$ solving the ODE $\deriv \Phi = \Phi\xi$.

The potential $\xi$ for the $\Xi_{g}^{n}$ families
is a Fuchsian potential on $\CPone$ with four simple poles
\begin{equation}
\xi \coloneq \sum_{k=0}^3 \frac{A_k\deriv z}{z-z_k}
\end{equation}
with order 2 symmetry
$\delta^\ast \xi = \sigma^{-1} \xi \sigma$,
$\delta(z) \coloneq -z$,
$\sigma \coloneq \diag(i,\,-i)$
and real symmetry
$\xi(\ol{z},\,\ol{\lambda}) = \ol{\xi}$.

The asymptotics of $\xi$ in the spectral parameter $\lambda$
is determined by the requirements of the generalized
Weierstrass representation and the Hopf differential of the
surface:
the upper right entry of the residue $A_0$ has a simple pole at
$\lambda=0$, the lower left entry of $A_0$ has a simple zero at
$\lambda=0$, and the residues of $\xi$ have no other poles
in the unit disk in $\bbC_\lambda$.

To reach the $\Xi_{g}^{n}$ families, we impose the condition
that the eigenvalues $\pm \nu_0$, $\pm \nu_1$,
be real and $\lambda$-independent, with
$\nu_0\in(0,\,1/4]$, $\nu_1\in[1/4,/1)$, $\nu_0 + \nu_1 = 1/2$.
This condition arises due to the fact that the eigenvalues control
the angle $\alpha$ of the rotational symmetry being opened.
With these assumptions, the
potential $\xi$ is a simpler replacement for the potential described
in equation~2.1 in~\cite{Heller_Schmitt_2015} , to which it is gauge
equivalent.

We note that the points in the punctured unit disk
in $\bbC_\lambda$ at which the parabolic structure corresponding
to $\xi$ is unstable are those points for which
the two eigenlines
of $A_1$ and $A_2$
corresponding to the positive eigenvalues
coincide.
We observed that the family $\Xi_{g}^{n}$, $g>1$
has $n$ unstable points in the unit disk.

\subsubsection{Geometric parameters}
The data determining a CMC surface in $\bbS^3$
via the generalized Weierstrass representation
is its potential $\xi$,
two sym points in $\bbS_\lambda^1$,
and the initial condition for the ODE $\deriv \Phi = \Phi\xi$.
For $\Xi_{g}^{n}$ families,
this data consists of
three geometric parameters
together with accessory parameters (coefficients of the residues of $\xi$).
The initial condition for the ODE $\deriv \Phi = \Phi\xi$
is determined
as the diagonal unitarizer of the monodromy,
unique up to isometry of $\bbS^3$.

The three real geometric parameters $(\gamma,\,\alpha,\,H)$ are
\begin{itemize}
\item
  the angle $\alpha \coloneq 4\pi\nu_0$ of the rotational symmetry being opened;
\item
  the conformal type $\gamma\coloneq [z_0,\,-z_1,\,-z_0,\,z_1]\in\bbR$ of the
  four punctured
  $\CPone$,
\item
  the mean curvature $H \coloneq i(\lambda_0+\lambda_1)/(\lambda_0-\lambda_1)$.
\end{itemize}

\subsubsection{Accessory parameters}
The eigenvalues of the residues of $\xi$,
determining the angle of the rotational symmetry,
must be controlled during the flow.
Since the complex dimension of the space of
monodromy representations on the 4-punctured sphere
with fixed eigenvalues is,
roughly speaking, $2$,
the residues of $\xi$ can be parametrized by
two meromorphic functions $\hat x$ and $\hat y$ of $\lambda$.
For numerical calculations, we represent
$\hat x$ and $\hat y$ as truncated power series
in $\lambda$ at $\lambda=0$.
Because the monodromy of $\xi$
is to be evaluated on the unit circle $\bbS_\lambda^1$,
we require that the potential $\xi$ is holomorphic
in the punctured unit disk. This holomorphicity
is achieved by the introduction of polynomials
and constraints on these polynomials.

More precisely, let $\hat x,\,\hat y$ be functions on the
unit disk with $\hat x$, $1/\hat x$ and $\hat y$ holomorphic.  Let
$p_k$, $q_k$, ($k\in\{0,\dots,3\}$) be polynomials satisfying the
constraints that $p_k$ monic and
\begin{equation}
  \label{eq:eigenvalue-constraint-0}
  e_{jk} \coloneq (p_j q_j - \nu_j) - (p_k q_k - \nu_k)
\spacecomma\quad (j,\,k\in\{0,\dots,3\})
\end{equation}
vanishes.
Under this constraint, the residues of the potential $\xi$ can be parametrized
by $p_k,\,q_k,\hat x,\,\hat y$, ($k\in\{0,\dots,3\}$) as
\begin{gather}
A_0 = \begin{bmatrix}-y & p_0p_2/\lambda\\-y_0 y_2\lambda & y\end{bmatrix}
  \spacecomma\quad
  A_1 = \begin{bmatrix}y & -y_1y_3/\hat{x}\\p_1 p_3\hat{x} & y\end{bmatrix}
    \\
p\coloneq p_0p_1p_2p_3\spacecomma\quad
y_k \coloneq q_k + p\hat y/p_k\spacecomma\quad
y \coloneq \nu_k + p_k y_k\spacecomma
\end{gather}
where $\nu_2 \coloneq -\nu_0$ and $\nu_3 \coloneq -\nu_1$.

For numerical computation, the accessory parameters are
\begin{equation}
A \coloneq (\operatorname{coeff} p_0,\dots\operatorname{coeff} p_3 \,|\,
\operatorname{coeff} q_0,\dots,\operatorname{coeff} q_3 \,|\,
\hat x_0,\dots,\hat x_N \,|\, \hat y_0,\dots,\hat y_N)
\end{equation}
where $\operatorname{coeff} q$ denotes the coefficients of a polynomial $q$,
and the series
\begin{equation}
\hat x = \sum_{k=0}^\infty \hat x_k\lambda^k
\spacecomma\quad
\hat y = \sum_{k=0}^\infty \hat y_k\lambda^k
\end{equation}
are truncated to power $N$.
The constraints~\eqref{eq:eigenvalue-constraint-0} are
\begin{equation}
  C_A \coloneq (\operatorname{coeff} e_{01},\,\operatorname{coeff} e_{02},\,
  \operatorname{coeff} e_{03})
\spaceperiod
\end{equation}

\subsection{The flow}
The generalized Whitham flow
is defined to preserve intrinsic and extrinsic closing conditions.
This flow is an implicit infinite dimensional ODE
computed numerically by truncation to a finite implicit system of the
form $A\dot X + B = 0$; $\dot X$ is obtained as the least squares
solution to this system.  The coefficients of the system depend on the
monodromy of the potential $\xi$, computed by a separate nested ODE.

\subsubsection{Intrinsic closing condition}
The intrinsic closing condition is that the monodromy of $\xi$ is
unitarizable along the unit circle $\bbS_\lambda^1$.  More concretely,
let $M_k$ ($k\in\{0,\dots,3\}$) be the monodromy generators for $\xi$
along a curve based at $z=0$ which winds once counterclockwise around
the pole $z_k$, and let $t_{jk} \coloneq \half\tr M_j M_k$.  By
proposition~2 in \cite{Heller_Schmitt_2015}, the monodromy
is unitarizable when $t_{ij}\in(-1,\,1)$.  Hence we impose the
constraint along $\bbS_\lambda^1$
\begin{equation}
  \label{eq:extrinsic-closing-condition-0}
  c_I \coloneq (\Imag t_{01}\spacecomma \Imag t_{02}\spacecomma \Imag t_{03}\spacecomma
  \Imag t_{12}\spacecomma \Imag t_{13}\spacecomma \Imag t_{23})
\spaceperiod
\end{equation}

For numerical computation, this constraint is implemented by imposing
the condition~\eqref{eq:extrinsic-closing-condition-0} at $S$
equidistant sample points $\mu_k=e^{2\pi i k/S}$, ($k\in\{0,\dots
S-1\}$) along $\bbS^1$.  For the flow, the number $S$ of sample points
must be large relative to the number $N$.  This is the vanishing of
\begin{equation}
C_I \coloneq (c_I(\mu_0),\,\dots,c_I(\mu_{S-1}))
\spaceperiod
\end{equation}

\subsubsection{Extrinsic closing condition}
The extrinsic closing conditions are that every monodromy of the
unitary frame $M$ satisfies $M(\lambda_0) = M(\lambda_1) \in \{\pm
\id\}$ at the two sympoints $\lambda_0,\,\lambda_1\in\bbS^1$.  By
proposition~1 in \cite{Heller_Schmitt_2015}, this is the condition
that
\begin{equation}
c_E \coloneq [\ell_0,\,\ell_1,\,\ell_2,\,\ell_3] - [z_0,\,-z_1,\,-z_0,\,z_1]
\end{equation}
vanishes to first order in $\lambda$ at each of the two sympoints,
where $\ell_k\in\CPone$ ($k\in\{0,\dots,3\}$)
is the eigenline of $A_k$ corresponding to its positive
eigenvalue.
With prime denoting the derivative with respect to $\lambda$,
this is the vanishing of
\begin{equation}
  C_{E} \coloneq \bigl( c_E(\lambda_0),\,c_E'(\lambda_0),\,
  c_E(\lambda_1),\,c_E'(\lambda_1) )
  \spaceperiod
\end{equation}

\subsubsection{The flow}
The flow is defined in terms of the function $f$, defined to vanish
when the intrinsic and extrinsic closing conditions are satisfied:
\begin{equation}
\begin{bmatrix}
  \text{geometric parameters $(\gamma,\,\alpha,\,H)$}\\
  \text{accessory parameters $A$}
\end{bmatrix}
\shortstack{$f$\\$\longmapsto$}
\begin{bmatrix}
  \text{intrinsic closing condition $C_I$}\\
  \text{extrinsic closing condition $C_E$}\\
  \text{accessory parameter constraint $C_A$}
\end{bmatrix}
\spaceperiod
\end{equation}
To flow along a curve in the 2-dimensional isosurface $f^{-1}(0)$
we consider
\begin{equation}
t\  \shortstack{$Y$\\$\longmapsto$}\ (t,\,u,\,A)
\ \shortstack{$h$\\$\longmapsto$}\ (\gamma,\,\alpha,\,H,\,A)
             \ \shortstack{$f$\\$\longmapsto$}\ (C_I,\,C_E,\,C_A)
\spacecomma
\end{equation}
where $t$ is the real flow parameter, $u$ is a real free parameter,
$h$ is an explicit immersion controlling the direction and speed of
the flow, and $Y$ is the sought function defined implicitly by the
condition $F\circ Y = 0$, where $F\coloneq f\circ h$.  With dot
denoting the derivative with respect to $t$, $Y$ is defined by the
implicit ODE
\begin{equation}
\deriv F \, \dot Y = 0
\spaceperiod
\end{equation}
In matrix form,
\begin{equation}
\begin{bmatrix}B & A\end{bmatrix}
  \begin{bmatrix}1\\\dot X\end{bmatrix} = 0
    \spacecomma
    \quad\text{that is}\quad
    A \dot X + B = 0
    \spaceperiod
\end{equation}
The vector field $\dot X$
is obtained as the least squares solution to the system
$A \dot X + B = 0$.

The generalized Whitham flow
starts with
the initial data for a homogeneous or Delaunay
torus~\cite{Heller_Heller_Schmitt_2015a},
with $(\gamma,\,\alpha,\,H) = (\text{constant},\,t,\,u)$,
reaching a surface in $\Xi_{g}^{n}$ by increasing genus.
Starting from such a surface,
the Whitham flow moves along the $\Xi_{g}^{n}$ family,
with $(\gamma,\,\alpha,\,H) = (t,\,\text{constant},\,u)$.

\subsection{Lawson symmetric surfaces}

A Lawson symmetric CMC surface is a compact CMC surface in $\bbS^3$ of
genus $g\ge 1$ which enjoys a cyclic symmetry of order $g+1$ with four
fixed points.  The families $\Xi_{g}^{n}$ described
above are Lawson symmetric,
and have an additional symmetry induced by the hyperelliptic
involution.  We conjecture,
for each pair of integers $g\ge 1$ and
$n\ge 1$, the existence of an additional 1-parameter family
$\widehat\Xi_{g}^{n}$ of Alexandrov embedded Lawson symmetric surfaces which
lack the symmetry induced by the hyperelliptic involution.  This
family is reachable via the generalized Whitham flow from the
$(2n+1)$-lobed Delaunay tori, and converges to a chain of $(g+1)n + 1$ CMC spheres.

\begin{conjecture*}
The space of Alexandrov embedded Lawson symmetric CMC surfaces
consists of the families $\Xi_{g}^{n}$ and $\widehat\Xi_{g}^{n}$.
\end{conjecture*}

In the case $g=1$, this moduli space is connected.
In the case of fixed $g>1$, the families $\Xi_{g}^{n}$ ranging over $n\in\bbN$
are disconnected from each other.


\end{document}